\documentclass[13pt]{article}
\usepackage{amsmath}
\usepackage{amscd}
\usepackage{amsfonts}
\usepackage{amssymb}
\usepackage{amsthm}
\usepackage{epsfig}
\usepackage{slashbox}
\usepackage[all]{xy}
 \newtheorem{lem}{Lemma}[section]
 \newtheorem{cor}[lem]{Corollary}
 
 \newtheorem{prop}[lem]{Proposition}
 \theoremstyle{definition}
 
 \theoremstyle{remark}

\begin{document}

\title{On the Signature Calculus for finite fields of order square of prime numbers}

\author{Qizhi Zhang \thanks{E-mail: zqz@ms.u-tokyo.ac.jp} \\}

\maketitle

\begin{abstract}
In \cite{Huang}, the authors proved that the discrete logarithm problem in a prime finite field is random polynomial time equivalent to computing the ramification signature of a real quadratic field. In this paper, we do this for a quadratic extension of a prime field. 
\end{abstract}



\section{Introduction}

Let $p$ be a big prime number and $E$ be an elliptic curve over $\mathbb{F}_p$. The computational complexity of the discrete logarithm problem in $\mathbb{F}_p ^\times$ or in $E(\mathbb{F}_p)$ is used in public-key encryption schemes. Estimating the greatest lower bound in the complexity of solving the discrete logarithm problem is an important and difficult problem.

In the case of $\mathbb{F}_p^\times$, the best algorithm known so far is the number field sieve. For example, see \cite{Gordon} and \cite{Schirokauer 1993}. It solves the discrete logarithm problem in $\mathbb{F}_p^\times$ in a conjectural running time $L_p(\frac{1}{3}, c) =\exp((c+
o(1))(\log p)^{1/3}(\log \log p)^{2/3}))$, with $c = (64/9)^{1/3}$. In the case of $E(\mathbb{F}_p)$, the problem can be reduced to the discrete logarithm problem in $\mathbb{F}_q^\times$ using the MOV attack (\cite{MOV}), where $q=p^k$ for some integer $k$. The discrete logarithm problem in $\mathbb{F}_q^\times$ can be solved  using the function field sieve ( \cite{Adleman 1994} ; \cite{Adleman and Huang 1999}) or a modified number field sieve ( \cite{Schirokauer 2000}). Let $e$ be the real number such that $k=(\log q / \log \log q)^e$. Then, the running time of the function field sieve is conjecturally equal to $L_q(\max\{ \frac{1}{3}, 1-e\}, O(1))$, and that of the modified number field sieve is conjecturally equal to $L_q(\max \{ \frac{1}{3}, \frac{1+e}{4}\}, O(1))$.

 Alternatively, we can estimate the greatest lower bound by studying an equivalent problem of a discrete logarithm problem. In \cite{Huang}, the authors lifted the discrete logarithm problem in $\mathbb{F}_p^\times$ to a real quadratic field. They defined the \lq\lq ramification signature" for the real quadratic field and proved that the discrete logarithm problem in $\mathbb{F}_p ^\times$  is random polynomial time equivalent to computing the ramification signature of the real quadratic field under two heuristic assumptions, namely, an assumption on the class number and an assumption on  a global unit of the real quadratic field.   

     In this paper, we lift the discrete logarithm problem in $\mathbb{F}_{p^2}^\times$ to a real quadratic field. We then define the \lq\lq ramification signature" for the real quadratic field and prove that the discrete logarithm problem in $\mathbb{F}_{p^2} ^\times$  is random polynomial time equivalent to computing the ramification signature of the real quadratic field, with one heuristic assumptions on the class number. We also show that in the proof of the equivalence in \cite{Huang} one can remove the assumption on the global unit. More precisely, we give an improvement ( Step 4 in section 3.2.b in the text ) on the construction of real quadratic field and global unit that makes the condition in  proposition 2 in section 4.1 in \cite{Huang} be satisfied automatically. 

In section 2, we define the ramification signature for a real quadratic field. In section 3, we prove the equivalence of the discrete logarithm problem in $\mathbb{F}_{q^2} ^\times$ and the computation of a ramification signature of a real quadratic field. Consequently, we also prove  the equivalence in \cite{Huang} without the assumption on the global unit.

\section{Signature}

To define the ramification signature for a real quadratic field, we need a proposition. 

\begin{prop}
Let $l$ and $p$ be two distinct odd prime numbers, and $K=\mathbb{Q} (\sqrt{D})$ be a real
quadratic field that splits over $l$ and inerts over $p$. We denote the ring of integers in $K$ by $A_K$, the points over $l$ by $u$ and $\tilde{u}$,  and the point over $p$ in $Spec A_K$ by $v$. Let $I_u, I_{\tilde{u}}$, and $I_v$ be the prime ideals of $A_K$ corresponding to $u, \tilde{u}$, and $v$, respectively. Let $Z:=\{ u, v \}$, and  $U:= Spec A_K \setminus Z$. Let $A_u$ and $A_v$ be the completions of $A_K$ at $u$ and $v$ respectively. Denote $A_U := \Gamma (U, \mathcal{O}_U)$ .

 Suppose that  $p^2-1$ is divisible by $l$, the class number of $K$ is not divisible by $l$, and  there is a unit $\alpha \in A_K ^\times$, such that 
\begin{displaymath}
\alpha ^{l-1} \neq 1 \mod I^2 _u , \quad \alpha ^{\frac{p^2-1}{l}} \neq 1 \mod I_v  .
\end{displaymath}  
Then, we have the following:

{\bf a.} There is an exact sequence 
 
$$
\xymatrix{
1 \ar[r] &  A_K^\times /A_K ^{\times l} \ar[r]^i & A_u^\times / A_u ^{\times l} \oplus A_v ^\times / A_v ^{\times l} \ar[r]^j & \pi_1(U)^{ab}/\pi _1(U)^{ab^l} \ar[r] &  1
}.            \eqno{(2.2.1)}
$$

{\bf b.} $\dim _{\mathbb{Z}/l\mathbb{Z}} \pi _1(U)^{ab}/\pi_1(U)^{ab^l}  =1$, where the $\pi_1(U)$ is the etale fundamental group of $U$ (see, for example, \cite{Etale});

{\bf c.} For any nontrivial character  $\chi: \pi_1(U) \longrightarrow \mathbb{Z}/l\mathbb{Z}$, $\chi$ is ramified at both $u$ and $v$.
\end{prop}

{\bf Proof.}

{\bf a.} Let us consider the following commutative diagram: 
\begin{displaymath}
\xymatrix{
1 \ar[r] & 1 \ar[r] \ar[d] & K^\times \ar[r]^\sim \ar[d] & K^\times \ar[r] \ar[d] & 1 \\
1 \ar[r] & \{ \pm 1 \} ^{\oplus 2} \oplus A_u^\times/A_u^{\times l} \oplus  A_v^\times/A_v^{\times l }\ar[r]  & \{ \pm 1\} ^{\oplus 2} \oplus K_u^\times/A_u^{\times l} \oplus K_v^\times/A_v^{\times l} \oplus \oplus_{x \neq u, v} \mathbb{Z} \ar[r]  & Div(K) \ar[r] & 0 .
}
\end{displaymath}
Through the snake lemma and class field theory, we have the following exact sequence:
\begin{displaymath}
\xymatrix{
A_K^\times \ar[r] & \{ \pm 1 \}^{\oplus 2} \oplus A_u^\times/A_u^{\times l} \oplus A_v^\times /A_v ^{\times l} \ar[r] & \pi _1(U)^{ab}/Im(A_u^{\times l} \oplus A_v^{\times l}) \ar[r] & Cl(K) \ar[r] & 1
} 
\end{displaymath} 
where the term $Im(A_u^{\times l} \oplus A_v^{\times l})$ is the image of $A_u^{\times l} \oplus A_v^{\times l}$ under the reciprocity map 
\begin{displaymath}
  K_u^\times \oplus K_v^{\times} \longrightarrow \mbox{Gal}(K_u^{ab}/K_u) \oplus \mbox{Gal}(K_v^{ab}/K_v) \longrightarrow  \pi_1(U)^{ab}.
\end{displaymath}

As the class number of $K$ is assumed  non-divisible by $l$, a diagram chasing shows an exact sequence

 \begin{displaymath}
A_K^\times /A_K ^{\times l} \longrightarrow A_u^\times / A_u ^{\times l} \oplus A_v ^\times / A_v ^{\times l} \longrightarrow \pi _1 (U)^{ab}/ {\pi _1 (U)^{ab}}^l \longrightarrow 1 .
\end{displaymath}
The hypothesis on the existence of the global unit shows that the left morphism is nonzero. Thus it is injective, since $A_K^\times /A_K ^{\times l}$ is an $\mathbb{Z}/l\mathbb{Z}$- linear space of dimension $1$. Therefore, we obtain the exact sequence $(2.2.1)$.

{\bf b.}  The complete discrete valuation rings $A_u$ and $A_v$ are isomorphic to $\mathbb{Z}_l$ and $W(\mathbb{F}_{p^2})$ ( the witt ring over $\mathbb{F}_{p^2}$) respectively. Therefore, the middle term in the sequence $(2. 2. 1)$ is isomorphic to $(\mathbb{Z}/l^2\mathbb{Z})^\times / (\mathbb{Z}/l^2\mathbb{Z})^{\times l} \oplus \mathbb{F}_{p^2}^\times / \mathbb{F}_{p^2}^{\times l}$, and is of $\mathbb{Z}/l\mathbb{Z}$-dimension $2$. Since we know that $A_K^\times/A_K^{\times l}$ is a $\mathbb{Z}/l\mathbb{Z}$-linear space of dimension 1, the right term in $(2. 2. 1)$ has $\mathbb{Z}/l\mathbb{Z}$-dimension 1. 

{\bf c.} We consider the dual sequence
$$
\xymatrix{
0 \ar[r] & Hom(\pi_1(U)^{ab}, \mathbb{Z}/l\mathbb{Z}) \ar[r]^{j^\star} &  Hom(A_u^\times / A_u ^{\times l} \oplus A_v ^\times / A_v ^{\times l} , \mathbb{Z}/l\mathbb{Z}) \ar[r] &  Hom(A_K^\times /A_K ^{\times l} , \mathbb{Z}/l\mathbb{Z}) \ar[r] & 0
}           \eqno{(2.2.2)}            
$$
of $(2.2.1)$. Denote the image of $\alpha$ under the morphism $i$ by $(\alpha _u, \alpha_v)$. For any $ \chi \neq 0 \in Hom(\pi_1(U)^{ab}, \mathbb{Z}/l\mathbb{Z})  $, denote the image of $\chi$ under the morphism $j^\star$ by $(\chi _u, \chi _v)$,  we then have
\begin{displaymath}
\langle \alpha _u, \chi_u \rangle+ \langle \alpha _v, \chi_v \rangle=0 \quad   \eqno{(2.2.3)}
\end{displaymath}
by $(2.2.2)$. Therefore, the following four conditions are equivalent:

(i).  $\chi$ is ramified at $u$,

(ii). $ \langle \alpha_u, \chi_u \rangle  \neq 0$,

(iii). $\langle \alpha _v, \chi_v \rangle \neq 0$,

(iv). $\chi$ is ramified at $v$.
 
The map $j^\star$ is injective, indicating that there is not non-trivial character $\chi: \pi_1(U) \longrightarrow \mathbb{Z}/l\mathbb{Z}$ such that it is unramified at both points $u$ and $v$. Therefore, for any non-trivial  $\chi \in Hom(\pi_1(U), \mathbb{Z}/l\mathbb{Z})$, $\chi$ must be ramified at both $u$ and $v$.   $\hfill \blacksquare$

The following corollary is proved in the proof of Proposition 2.1c.

\begin{cor}
Under the conditions in proposition 2.1, for any non-trivial character $\chi: \pi_1(U) \longrightarrow \mathbb{Z}/l\mathbb{Z}$, we have the following:

(i) $\langle \alpha _u, \chi_u \rangle \neq 0$,

(ii) $\langle \alpha _v, \chi_v \rangle \neq 0$,

(iii) $\langle \alpha _u, \chi_u \rangle + \langle \alpha _v, \chi_v \rangle =0$.

\end{cor}

Through the natural isomorphism $A_u^\times /A_u^{\times l} \cong (\mathbb{Z}/l^2\mathbb{Z})^\times / (\mathbb{Z}/l^2\mathbb{Z})^{\times l}$, $A_u^\times / A_u ^{\times l}$ is generated by $1+l$. For any generator $g$ of $\mathbb{F}_{p^2} ^\times / \mathbb{F} _{p^2} ^{\times l}$, we regard it as a generator of $A_v ^\times / A_v ^{\times l}$ through the natural isomorphism $A_v ^\times /A_v ^{\times l} \cong \mathbb{F}_{p^2}^\times / \mathbb{F}_{p^2}^{\times l}$. Clearly, $\langle 1+l, \chi _u \rangle ^{-1}\langle g, \chi_v \rangle$ is independent of the choice of $\chi \neq 0 \in Hom(\pi_1(U), \mathbb{Z}/l\mathbb{Z})$. We call this term the {\bf ramification signature of $U$ with respect to $g$}.

\section{Signature computation problem  and discrete logarithm problem in $\mathbb{F}_{p^2}^\times $}

In this section, we show that the discrete logarithm problem in $\mathbb{F}_{p^2}^\times$ is random polynomial time equivalent to computing the ramification signature of some real quadratic field.

\subsection{Reduction from signature computation problem to discrete logarithm problem}

Suppose given $p, l$, $K=\mathbb{Q} (\sqrt{D})$, $U$, $u, \tilde{u}, v$, $\alpha$, $g$, as in Proposition 2.1. Then the computation of the ramification signature of $U$ with respect to $g$ can be reduced to a discrete logarithm problem in $\mathbb{F}_{p^2}$ as follows by using Corollary 2.2.

Let us consider the following commutative diagram: 
\begin{displaymath}
\xymatrix{
A_K ^\times \ar[r] & A_u ^\times \ar[r]^\sim & \mathbb{Z}_l ^\times \ar[r] \ar[d]  & \mathbb{Z}_l^\times / \mathbb{Z}_l ^{\times l} \ar[d] \\
& & (\mathbb{Z}/l^2 \mathbb{Z})^\times \ar[r] &   (\mathbb{Z}/l^2\mathbb{Z})^\times / (\mathbb{Z}/l^2\mathbb{Z})^{\times l}
.}
\end{displaymath}   
If the image in $(\mathbb{Z}/l^2 \mathbb{Z})^\times$ of $\alpha$ equals $\xi (1+l)^y$, where $\xi$ is an($l-1$)-st root of unity, then its image in $(\mathbb{Z}/l^2 \mathbb{Z})^\times/ (\mathbb{Z}/l^2 \mathbb{Z})^{\times l} $ will be equal to $(1+l)^y$. We can easily compute $\xi$, $y$ and consequently the first term in (2.2.3) $\langle \alpha _u, \chi_u \rangle =y \langle 1+l, \chi_v \rangle$.  

For the second term in (2.2.3), if the image of $\alpha$ under the morphism $A_K^\times \longrightarrow A_v ^\times /A_v ^{\times l} \cong \mathbb{F}_{p^2}^\times / \mathbb{F}_{p^2}^{\times l}$ is $a=g^m$, then $\langle \alpha_v, \chi_v \rangle =m \langle g, \chi_v \rangle $.

By Corollary 2.2, if we can compute $m$ from $a=g^m$, then
we can compute
\begin{displaymath}
\langle 1+l, \chi_u\rangle ^{-1}\langle g, \chi _v\rangle = -m^{-1}y  \quad \in \mathbb{Z}/l\mathbb{Z}  .
\end{displaymath}

\subsection{Reduction from discrete logarithm problem to signature computation problem }

Let  $g$ be a generator of $\mathbb{F}_{p^2}^\times$, $a \in \mathbb{F}_{p^2}^\times$ and $l$ be a prime dividing  $p^2-1$. The computation of discrete logarithm $\log _g a \mod l$ can then be reduced to  computing the ramification signature of a real quadratic field as follows, by using Corollary 2.2.

Let $m \equiv \log _g a \mod l$. If $a \in \mathbb{F}_{p^2}^{\times l}$, then we have $m \equiv 0 \mod l$. Thus we can suppose $a \notin \mathbb{F}_{p^2}^{\times l}$.

{\bf a.} If $l \nmid p-1$, we must have $l | p+1$. Let $\tilde{a}:=a^{p-1}, \tilde{g}:=g^{p-1}$. We then have
\begin{displaymath}
\tilde{a} = \tilde{g}^m, \quad Nm(\tilde{a})=Nm(\tilde{g})=1, \quad \tilde{a} \notin \mathbb{F}_{p^2}^{\times l} .
\end{displaymath}
We take $t \in \mathbb{F}_p$ such that  $(\frac{t}{p})=-1$. We have  $\mathbb{F}_{p^2}=\mathbb{F}(\sqrt{t})$. We put $\tilde{a}=a_0+b_0 \sqrt{t}$, where $a_0, b_0 \in\mathbb{F}_p$. We can assume $b_0 \neq 0$; otherwise, $\tilde{a}^2=1$ and $m=\frac{p+1}{2}$ or $p+1$.

We have $a_0^2-b_0 ^2t=Nm(\tilde{a})\equiv 1 \mod p$. Hence, for any $k \in \mathbb{Z}$, the following holds: 
\begin{displaymath}
(\frac{(a_0+kp)^2-1}{p})=(\frac{a_0^2-1}{p})=(\frac{b_0^2t}{p})=(\frac{t}{p})=-1 .
\end{displaymath}
We choose $k \in \{0, 1, \cdots , l-1 \}$ randomly, until $(\frac{(a_0+kp)^2-1}{l})=1$. Lemma 3.1 below for $c=1$  shows that we can obtain such $k$ with probability about 50\% each time.  

If we find such $k$, let $a_1:=a_0 +kp \in \mathbb{Z}^\times$. We have $\sqrt{a_1^2-1} \in \mathbb{Z}_l$ because $(\frac{a_1^2-1}{l})=1$. If $(a_1+\sqrt{a_1^2-1})^{l-1} \not \equiv 1 \mod l^2$, we know also $(a_1-\sqrt{a_1^2-1})^{l-1} \not \equiv 1 \mod l^2$, let $x=a_1$. Else, let $x=a_1+pl$. Lemma 3.2 below for $c=1$  shows that $(x+\sqrt{x^2-1})^{l-1} \not \equiv 1 \mod l^2$ and $(x-\sqrt{x^2-1})^{l-1} \not \equiv 1 \mod l^2$ .

Let $K:=\mathbb{Q}(\sqrt{x^2-1})$. Then,  $K$ inerts over $p$ and splits over $l$  because $(\frac{x^2-1}{p})=-1$, and $(\frac{x^2-1}{l})=1$. Let $v \in Spec A_K$ be the point over $p$ and $u \in SpecA_K$ be a point over $l$. We have $\sqrt{x^2-1}\equiv \pm b_0 \sqrt{t} \quad \mod v$, because $x^2-1 \equiv a^2-1 \equiv b_0^2t \quad \mod v$. 
Let $\alpha:= x+\sqrt{x^2-1} \in A_K$, if $\sqrt{x^2-1}\equiv \pm b_0 \sqrt{t} \quad \mod v$; or $\alpha:= x-\sqrt{x^2-1} \in A_K$, if else. We then have $\alpha ^{l-1} \not \equiv 1 \mod I_u^2$  and 
\begin{displaymath}
\alpha \equiv a_0+\sqrt{a_0^2-1} \equiv a_0+b_0 \sqrt{t} \equiv \tilde{a} \equiv \tilde{g}^m \quad \mod v
\end{displaymath}
implying that  $\alpha ^\frac{p^2-1}{l} \not \equiv 1 \mod I_v$ as  $\tilde{a} \notin \mathbb{F}_{p^2}^{\times l}$. As  $\alpha:= x+\sqrt{x^2-1} \in A_K$ and $Nm(\alpha)=x^2-(x^2-1)=1$, we have $\alpha \in A_K^\times$.

Let $U:= Spec A_K \setminus \{u, v\}$. We assume  that $l \nmid h_K$, which is likely  to be  satisfied. Proposition 2.1 then shows $\langle \alpha _u, \chi \rangle + \langle \alpha_v, \chi \rangle =0$, for any $\chi \neq 0 \in Hom(\pi_1(U), \mathbb{Z}/l\mathbb{Z})$. Let $(1+l)^y$ be the image of $\alpha$ under the morphism
\begin{displaymath}
A_K^\times \longrightarrow A_u^\times \cong \mathbb{Z}_l^\times \longrightarrow (\mathbb{Z}/l^2\mathbb{Z})^\times /(\mathbb{Z}/l^2\mathbb{Z})^{\times l}.
\end{displaymath}
For the first term in $(2.2.3)$, we have  $\langle \alpha _u, \chi \rangle=y\langle 1+l, \chi \rangle$. For the second term in $(2.2.3)$, we have $\langle \alpha _v, \chi \rangle=m \langle \tilde{g}, \chi \rangle$. By Corollary 2.2.(iii), we obtain
\begin{displaymath}
y\langle 1+l, \chi \rangle +m \langle \tilde{g}, \chi \rangle=0  .
\end{displaymath}
Therefore, if we can compute the ramification signature $\langle \chi, 1+l \rangle ^{-1} \langle \chi, g \rangle $ of $U$ with respect to $g$, then we can compute $m=-y \langle \chi, 1+l \rangle ^{-1}\langle \chi, g \rangle$.

{\bf b.} If $l | p-1$, we have $Nm(a)=Nm(g)^m$ as elements in $\mathbb{F}_p$. The construction in \cite{Huang} gives us a real quadratic field and a global unit in the field that enable us to reduce the computation of $m$ satisfying $Nm(a)=Nm(g)^m$ to the signature computation problem of the real quadratic field using the algorithm in \cite{Huang}. However, the construction requires some conditions on the class number of the field and the unit to be satisfied (\cite{Huang}; Section 4.2). We give an improvement in Step 4 below on the construction recalled below. With the improvement, one of the condition (the condition 2 in Proposition 2 of section 4.1 in \cite{Huang}) is satisfied automatically.

Let $\tilde{a}=\tilde{g}^m$ in $\mathbb{F}_p^\times$ where $m$ is to be computed. If $\tilde{a}^\frac{p-1}{l}=1$, then $m \equiv 0 \quad \left(\mod l \right) $.  Thus suppose $\tilde{a}^{\frac{p-1}{l}} \neq 1$. We will lift $\tilde{a}$ to some unit $\alpha$ of a real quadratic field $K$ such that $\alpha \equiv \tilde{a}  \mod v $ for some place $v$ of $K$ over $p$, $\alpha ^{l-1} \neq 1 \mod I_u^2$, and $\alpha ^{l-1} \neq 1 \mod I_{u'}^2$ for the two places $u$ and $u'$ of $K$ over $l$. We do it as follows: 

1. Compute $\tilde{b} \in \mathbb{F}_p^\times$ such that $\tilde{a}\tilde{b}=1$ in $\mathbb{F}_p^\times$. 

2. Put  $c:=\frac{\tilde{a}+\tilde{b}}{2}$, $d:=\frac{\tilde{a}-\tilde{b}}{2}$. Note that $c^2-d^2=1$ and $\tilde{a}=c+d$. We can assume $d \neq 0$; otherwise, $\tilde{a}^2=1$ and $m=\frac{p-1}{2}$ or $p-1$.

3. Lift $d$ to an integer. We have $(\frac{1+d^2}{p})=(\frac{c^2}{p})=1$. We choose $k \in \{ 0, 1, ..., l-1 \}$ randomly until $(\frac{(d+kp)^2+1}{l})=1$. Lemma 3.1 below for $c=-1$ shows that we can obtain such $k$ with probability of  about 50\% each time. 

4. If we find such $k$, let $d_1:=d +kp \in \mathbb{Z}_l^\times$. We may take $\sqrt{d_1^2+1} \in \mathbb{Z}_l^\times$ since $(\frac{d_1^2+1}{l})=1$. If $(d_1+\sqrt{d_1^2+1})^{l-1} \equiv 1 \mod l^2$, let $x=d_1$; otherwise let $x=d_1+pl$. Lemma 3.2 below  for $c=-1$  shows that $(x+\sqrt{x^2+1})^{l-1} \not \equiv 1 \mod l^2$,  $(x-\sqrt{x^2+1})^{l-1} \not \equiv 1 \mod l^2$.

5. Let $K:=\mathbb{Q}(\sqrt{x^2+1})$, $\alpha :=x+ \sqrt{x^2+1}\in \mathcal{O}_K$. Note that $Nm(\alpha)=1$, so $\alpha$ is a unit of $K$.

6. Let  $v$ be the point in $Spec \mathcal{O}_K$ responding to the prime ideal $(p, \sqrt{x^2-1}-c)$, $v'$ be the point in $Spec \mathcal{O}_K$ responding to the prime ideal $(p, \sqrt{x^2-1}+c)$, $u$ and $u'$ be the points in $Spec \mathcal{O}_K$ over $l$. Thus, $\alpha \equiv d+c \equiv \tilde{a} \quad \mod v$, $\alpha \equiv d-c \equiv -\tilde{b} \quad  \mod v'$, $\alpha ^{l-1} \neq 1 \mod I_u^2 $ and $\alpha ^{l-1} \neq 1 \mod I_{u'}^2$    $\hfill \blacksquare$

In \cite{Huang}, they proved the reduction from a signature computation to a discrete logarithm problem in $\mathbb{F}_p$ without any heuristic assumption. Therefore,  we conclude that the  discrete logarithm problems in $\mathbb{F}_{p^2}$ and  $\mathbb{F}_p$ are  random polynomial time equivalent to  some signature computation problem with only one assumption, namely, that on the  class number.

The following are the statements and proofs of lemma 3.1 and lemma 3.2.

\begin{lem}
Let $l$ be an odd prime, $c \in \mathbb{F}_l^\times$. Define a map $f: \mathbb{Z}/l\mathbb{Z} \longrightarrow \{ 0,1,-1 \}$ by $a \mapsto (\frac{a^2-c}{l})$. Then, we have
\begin{displaymath}
|f^{-1}(0)|=2, \quad |f^{-1}(1)|=(l-3)/2,  \quad  |f^{-1}(-1)|=(l-1)/2  \quad \mbox{ if } (\frac{c}{l})=1,
\end{displaymath} 
\begin{displaymath}
|f^{-1}(0)|=0, \quad |f^{-1}(1)|=(l-1)/2,  \quad  |f^{-1}(-1)|=(l+1)/2  \quad \mbox{ if } (\frac{c}{l})=-1.
\end{displaymath}
\end{lem}
 
{\bf Proof.} Let $X$ be the curve defined by $y^2=x^2-c$ over $\mathbb{F}_l$. For any $a \in \mathbb{F}_l$, the cardinality of the set 
\{ $\mathbb{F}_l$-rational point of $X$ that has first coordinate $a$ \} is $f(a)+1$. Therefore, the following holds:
\begin{displaymath}
\sum _{a \in \mathbb{F}_l} (f(a)+1) = |X(\mathbb{F}_l) |.
\end{displaymath}
The curve $X$ is isomorphic to the affine scheme defined by $z\omega=1$ over $\mathbb{F}_l$, which implies $|X(\mathbb{F}_l)|=l-1$, and $\sum _{a \in \mathbb{F}_l} f(a) = |X(\mathbb{F}_l) |-l=-1$. Clearly,
\begin{displaymath}
\begin{array}{ll}
f^{-1}(0)=\{ \sqrt{c}, -\sqrt{c} \} &  \mbox{ if } (\frac{c}{l}=1),   \\
f^{-1}(0)=\phi & \mbox{ if } (\frac{c}{l}=-1).
\end{array} 
\end{displaymath}
and the lemma follows easily from $|f^{-1}{1}|+|f^{-1}(-1)|+|f^{-1}(0)|=l$ and $|f^{-1}(1)|-|f^{-1}(-1)|=\sum _{a \in \mathbb{F}_l} f(a)=-1$.     $\hfill \blacksquare$

\begin{lem}
Let $p$ and $l$ be two distinct odd prime numbers. Let $c$ be an integer such that $c^{l-1} \equiv 1 \mod l^2$  and  $a$ be an integer such that $(\frac{a^2-c}{l})=1$. We denote a square root of $a^2-c$ in $\mathbb{Z}_l^\times$ by $\sqrt{a^2-c}$. If $(a+\sqrt{a^2-c})^{l-1} \in 1+l^2\mathbb{Z}_l$, then we have $((a+pl)+\sqrt{(a+pl)^2-c})^{l-1} \notin 1+l^2\mathbb{Z}_l$ and  $((a+pl)-\sqrt{(a+pl)^2-c})^{l-1} \notin 1+l^2\mathbb{Z}_l$. 
\end{lem}

{\bf Proof.} 
By Hensel's lemma, there is a unique square root $\sqrt{(a+x)^2-c}$ of $(a+x)^2-c$ in $\mathbb{Z}_l[[x]]$ such that it's image under the morphism $x \mapsto 0: \mathbb{Z}_l[[x]] \longrightarrow \mathbb{Z}_l$ is $\sqrt{a^2-c}$. Let $h(x):=(a+x)+\sqrt{(a+x)^2-c}$, we then have
$$
h(pl) \equiv h(0) + h'(0)pl \quad \mod l^2 ,
$$
where $h'(0)=1+\frac{a}{\sqrt{a^2-c}}=\frac{h(0)}{\sqrt{a^2-c}}$. Therefore, we have
$$
h(pl) \equiv h(0)(1+\frac{p}{\sqrt{a^2-c}} l) \mod l^2    .
$$
The term $\frac{p}{\sqrt{a^2-c}}$ is not divided by $l$, which implies $h(pl)^{l-1} \not \equiv h(0)^{l-1} \mod l^2$. Hence, we have
\begin{displaymath}
\begin{array}{rl}
&((a+pl)+\sqrt{(a+pl)^2-c})^{l-1}  \\
\not \equiv & (a+ \sqrt{a^2-c})^{l-1} \quad \quad \mod l^2  \\ 
\equiv & 1  \quad \quad \mod l^2.    
\end{array}
\end{displaymath}

The fact that  $((a+pl)+\sqrt{(a+pl)^2-c})^{l-1} ((a+pl)-\sqrt{(a+pl)^2-c})^{l-1}=c^{l-1} \equiv 1 \mod l^2$ shows 

 $((a+pl)-\sqrt{(a+pl)^2-c})^{l-1} \not \equiv 1 \mod l^2$. $\hfill \blacksquare$

\noindent $\mathbf{Acknowledgments.}$ I would like to thank professor Takeshi
Saito for giving me valuable advice.


\end{document}